\documentclass[11pt]{revtex4}

\usepackage{graphicx,epsfig,amssymb,amsbsy,amsmath,amsthm}
\newcommand{\bb}[1]{\left({#1}\right)}					% ( )
\newcommand{\sq}[1]{\left[#1\right]}						% [ ]
\newcommand{\cc}[1]{\left\{#1\right\}}					% { }
\newcommand{\op}[1]{\mathcal{#1}}
\newcommand{\ord}[1]{{\sf O}\bb{#1}}					% order
\newcommand{\sfrac}[2]{\mbox{$\frac{#1}{#2}$}}	

\newcommand{\hf}{\mbox{$\frac12$}}

\newtheorem{theorem}{Theorem}

\newtheorem{proposition}[theorem]{Proposition}

\newcommand{\sign}{\operatorname{sign}}

\newcommand{\fref}[1]{figure~\ref{#1}}
\newcommand{\Fref}[1]{Figure~\ref{#1}}

\newcommand{\eref}[1]{(\ref{#1})}

\newcommand{\sref}[1]{section~\ref{#1}}

\def\eps{\varepsilon}
\begin{document}

\title{Equivalence between the nonsmooth ``two-fold'' and slow-fast ``folded'' flow singularies}
\title{Fold singularities of nonsmooth and slow-fast dynamical systems\\ -- equivalence by the hidden dynamics approach}
\author{Mike R. Jeffrey}\address{Dept. of Engineering Mathematics, University of Bristol, Merchant Venturer's Building, Bristol BS8 1UB, UK, email: mike.jeffrey@bristol.ac.uk}
%\affiliation{Department of Engineering Mathematics, University of Bristol, Merchant Venturer's Building, Bristol BS8 1UB, UK}
\date{\today}

\begin{abstract} 
The {\it two-fold singularity} has played a significant role in our understanding of uniqueness and stability in piecewise smooth dynamical systems. When a vector field is discontinuous at some hypersurface, it can become tangent to that surface from one side or the other, and tangency from both sides creates a two-fold singularity. The flow this creates bears a superficial resemblance to so-called {\it folded singularities} in (smooth) slow-fast systems, which arise at the intersection of attractive and repelling branches of slow invariant manifolds, important in the local study of canards and mixed mode oscillations. Here we show that these two singularities are intimately related. When the discontinuity in a piecewise smooth system is blown up or smoothed out at a two-fold singularity, the resulting system can be mapped onto a folded singularity. The result is not obvious, however, since it requires the presence of nonlinear or `hidden' terms at the discontinuity, which turn out to be necessary for structural stability of the blow up (or smoothing), and necessary for mapping to the folded singularity. 
\end{abstract}
%\pacs{}

\maketitle

\section{Introduction}\label{sec:intro}

For a local singularity that is so easy to define, the `two-fold' singularity of piecewise-smooth dynamical systems has proven surprisingly difficult to characterize. It has been the subject of interest both for its intricate phase portrait \cite{f88,jc12}, its ambiguous stability \cite{t90,jc09}, and recently because of its role in determinacy-breaking \cite{cj10}. Now that it seems these are all well understood, our interest in this paper is in how the local dynamics relates to that of smooth flows. %, in a rather more direct way than . 
%Consider a system of the form
%\begin{equation}\label{gen}
%{\bf \dot x}={\bf f}({\bf x};\lambda)\qquad\mbox{ where} \quad\lambda={\rm sign}\bb{h({\bf x})}
%\end{equation}
%${\bf \dot x}={\bf f}({\bf x};\lambda)$ where $\lambda={\rm sign}\bb{h({\bf x})}$, 
%where $\bf f$ is a vector field and $h$ a scalar function, so a discontinuity takes place in $\bf f$ along a smooth hypersurface $h=0$. 

\begin{figure}[h!]\centering\includegraphics[width=\textwidth]{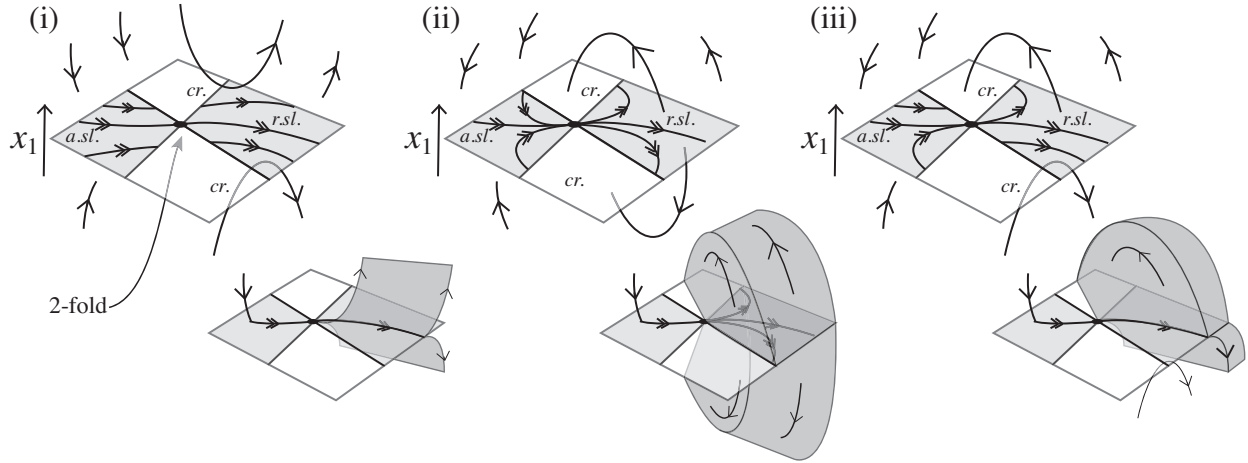}
\vspace{-0.3cm}\caption{\small\sf Three kinds of two-fold. The main figure shows the phase portrait: in the unshaded regions the flow crosses ($cr.$) through a discontinuity at $x_1=0$, in the shaded regions the flow can only {\it slide} along the discontinuity on $x_1=0$, the region being attracting ($a.sl.$) or repelling ($r.sl.$). In the examples shown, determinacy-breaking occurs at the singularity, meaning that the flow there becomes set-valued, the set has 2 dimensions in (i) and 3 dimensions in (ii-iii). }\label{fig:2folddet}\end{figure}

%Quite simply, i. 
If a discontinuity occurs in the vector field of a flow on some hypersurface $\Sigma$, and the flow curves (or `folds' parabolically) towards or away from $\Sigma$ on both sides of the surface, then under generic conditions the result is a two-fold singularity, as depicted in \fref{fig:2folddet}. 
%
%The dynamics around the two-fold is now understood in considerable detail. 

Starting with \cite{f88,t90}, the elucidation of the two-fold singularity's dynamics can be traced through references in \cite{jc12} (where the local form of the two-fold singularity is also extended to more than three dimensions). Though it can create limit cycles via a local bifurcation \cite{cj10}, the two-fold singularity itself is not an attractor. The main interest therefore lies in its transitory effects on a flow, particularly in cases where it forms a bridge {\it from} an attracting region on the discontinuity surface {\it into} a repelling region, mimicking {\it canard} behaviour of smooth two-timescale systems \cite{b81} but with a more extreme character, because the two-fold breaks determinacy in both forward and backward time through the singularity. This is illustrated in \fref{fig:2folddet} by a typical single trajectory entering the singularity, being deterministic until it does so, and afterwards exploding into a set-valued flow of infinite onward trajectories whose local form, however, is still somewhat constrained by the local geometry. 

Particularly because of some similarity to canard dynamics, attention has turned to how the two-fold can be understood as a limit or approximation of a smooth flow. The topological equivalence between ``sliding'' motion along a discontinuity surface, and ``slow'' motion on invariant manifolds of smooth two-timescale systems, has been shown \cite{ts11}. A qualitative connection has also been made between the sliding phase portraits at two-fold singularities and the slow manifold dynamics at so-called `folded' singularities of two-timescale systems \cite{jd10,jd12}. We shall prove here a more direct connection between the two singularities, by showing equivalence under explicit coordinate transformations.

The two singularities that concern us can be defined as isolated points in $n\ge3$ dimensions satisfying the following systems of equations:
\begin{equation}\label{sys}
\begin{array}{rllll}
&&\mbox{The two-fold singularity:}&&\mbox{The folded singularity:}\\
\rm system:&&\left\{\begin{array}{l}\dot x_i=f_i(x_1,x_2,...;\lambda)\\{\rm for}\;\;i=1,...,n\end{array}\right.&\qquad\qquad&
\left\{\begin{array}{l}\eps\dot y_1=g_1(y_1,y_2,...;\eps)\\\;\;\dot y_j=g_j(y_1,y_2,...;\eps)\\{\rm for}\;\;j=2,...,n\end{array}\right.\\
\rm parameter:&&\quad\lambda={\rm sign}\bb{x_1}&&\quad0<\eps\ll1\\
\rm singularity:&&\quad0=x_1=\left.\displaystyle\lim_{\delta\rightarrow0}\dot x_1\right|_{x_1=\pm\delta}&&\quad0=g_1=\dot g_1=\frac{\partial g_1}{\partial g_1}\;
\end{array}
\end{equation}
along with a number of non-degeneracy conditions, which state that immediate higher derivatives have typical (i.e. non-vanishing) values; these will be elaborated on later. The functions $f_i$ and $g_i$ are assumed to be differentiable with respect to their arguments. The parameter $\lambda$ introduces a discontinuity via the sign function, while $\eps$ introduces a separation into slow and fast timescales $t$ and $t/\eps$ (where $\dot z$ of course denotes the time derivative $\frac{d\;}{dt}z$). Three dimensions are sufficient to understand the basic dynamics of either singularity, so henceforth we restrict to $n=3$. 

In an attempt to relate the two systems in \eref{sys}, we might first observe that the folded singularity system becomes discontinuous in the singular $\eps\rightarrow0$ limit. However, that limit is {\it not} equivalent to the discontinuous system of the two-fold singularity in \eref{sys}. To understand how the two systems in \eref{sys} {\it are} related we must instead start with the two-fold system, and somehow regularize the discontinuity to understand what forms of smooth system have it as a limit. 

Regularizing or smoothing a discontinuity raises immediate issues of uniqueness, namely that infinitely many qualitatively different smooth systems can have the same discontinuous limit. We shall find that actually this ambiguity can be encoded into the discontinuous system by means of {\it nonlinear switching} terms --- nonlinear dependence on the parameter $\lambda$ in \eref{sys} --- as introduced in \cite{j13error} (sometimes called `hidden' terms because they vanish everywhere except at the discontinuity). We shall see that the smooth system obtained from a two-fold singularity is structurally unstable if it depends only linearly on $\lambda$, and we shall show that a small perturbation, by terms that are nonlinear in $\lambda$, restores structural stability and allows transformation into the general local singularity expected in a smooth system, namely the folded singularity. 

In \sref{sec:2fold} we introduce the normal form of the two-fold singularity, and outline the basic steps for its study by blowing up the discontinuity in \sref{sec:blowup}. 
In \sref{sec:2fdummy} we blow up the normal piecewise smooth system, assuming only linear dependence on $\lambda$ as in the standard literature, giving a system which we find is structurally unstable. 
In \sref{sec:2fpert} we introduce a nonlinear perturbation of the normal form, then blow this up to find that it removes the structurally instability, and can be mapped onto the folded singularity of a smooth two timescale system. Remarks showing that these results follow also if we smooth, rather than blow up, the discontinuity, are given in \sref{sec:conc}.

%A simple question can be asked about the many forms of dynamics that this singularity exhibits, and that is whether there is any analogue of the two-fold singularity in smooth systems. 

%We begin in \sref{sec:smooth} by showing the role that the nonlinear term $\bf g$ plays in one type of regularization --- smoothing. We then study determinacy breaking at a two-fold singularity in \sref{sec:2fold}, including how to resolve for dynamics at the discontinuity, and how to show that certain piecewise-smooth behaviours are stable to perturbation. 
%Some closing remarks are given in \sref{sec:conc}. 

%%%%%%%%%%%%%%%%% %%%%%%%%%%%%%%%%%%
%%%%%%%%%%%%%%%%% %%%%%%%%%%%%%%%%%%
%%%%%%%%%%%%%%%%% %%%%%%%%%%%%%%%%%%
\section{The two-fold singularity}\label{sec:2fold}

The normal form of the two-fold singularity is 
\begin{equation}\label{2f}
\bb{\dot x_1,\dot x_2,\dot x_3}=%\bb{f_1,f_2,f_3}\;=\;
\left\{\begin{array}{lll}\bb{-x_2,a_1,b_1}&\rm if&x_1>0\;,\\\bb{+x_3,b_2,a_2}&\rm if&x_1<0\;,\end{array}\right.
\end{equation}
in terms of constants $a_i=\pm1$ and $b_i\in\mathbb R$. By results in \cite{f88,t93,jc12}, a system is locally approximated by \eref{2f} when it satisfies the conditions in the lefthand column of \eref{sys}, along with non-degeneracy conditions stating that the flow curvature is quadratic, given by $\displaystyle\lim_{\delta\rightarrow0}\left.\ddot x_1\right|_{x_1=\pm\delta}\neq0$, and that the flows either side of the discontinuity are independent of each other and of the discontinuity surface, meaning the vectors along the $x_1$ direction, along $\displaystyle\lim_{\delta\rightarrow0}\left.\dot x_1\right|_{x_1=+\delta}$, and along $\displaystyle\lim_{\delta\rightarrow0}\left.\dot x_1\right|_{x_1=-\delta}$, are linearly independent. 

The local flow `folds' towards or away from the switching surface $x_1=0$, along the line $x_2=x_1=0$ on one side of the surface, and along the line $x_3=x_1=0$ on the other. Hence the point where these lines cross is called the `{\it two}-fold'. As a result, the surface $x_1=0$ is attractive in $x_2,x_3>0$ and repulsive in $x_2,x_3<0$, while trajectories cross the surface transversely in $x_2x_3<0$. In the attractive and repulsive regions the flow {\it slides} along the surface $x_1=0$, and the vector field it follows is found by interpolating across the discontinuity; we shall see how this is done below. 

The qualitative picture is then as shown in \fref{fig:2folddet}. 
%For various classes of the parameters $a_{1,2}$ and $b_{1,2}$ this system exhibits determinacy-breaking as depicted in \fref{fig:2folds}. 
The precise form of the local dynamics depends on whether the flow curves towards or away from the discontinuity, determined by $a_1$ and $a_2$, and also depends crucially on the quantity $b_1b_2$, which quantifies the jump  in the angle of the flow across the discontinuity. An accounting of the many classes of dynamics that arise from these simple conditions is given in \cite{jc12} and references therein, we give only the pertinent details here.

%The phase portrait that surrounds a two-fold singularity takes a number of different forms depending on the signs of $a_1,a_2,$ and the values of $b_1,b_2$%; they can be found in \cite{jc12}, so we sketch only some of the more interesting behaviours in \fref{fig:2folddet}
%. 
The three main `flavours' of two-fold are: the visible two-fold for $a_1=a_2=-1$, the invisible two-fold for $a_1=a_2=1$, and the mixed two-fold for $a_1a_2=-1$; an example of each is shown in \fref{fig:2folddet} (i,ii,iii) respectively. In the cases depicted, there exist one or more trajectories passing from the attractive sliding region to the repelling sliding region, called {\it canard} orbits. This passage occurs in finite time (since the vector field \eref{2f} is non-vanishing everywhere locally). The flow is unique in forward time everywhere except in the repelling sliding region, where it is set-valued because trajectories may slide along $x_1=0$, but also be ejected into $x_1>0$ or $x_1<0$ at any point. This means that the flow may evolve deterministically until it arrives at the singularity by means of a canard, at which point it becomes set-valued, so we say that {\it determinacy breaking} occurs at the singularity whenever canards are exhibited. This occurs in the invisible case when $b_1,b_2<0$ and $b_1b_2>1$, in the visible case when $b_1<0$ or $b_2<0$ or $b_1b_2<1$, and finally in the mixed case when $b_1<0<b_2$ and $b_1b_2<-1$ or when $b_1+b_2<0$ and $b_1-b_2<-2$.
(The particular cases shown in \fref{fig:2folddet} are: (i) $a_1=a_2=-1$ with $b_1<0$ or $b_2<0$ or $b_1b_2<1$; (ii) $a_1=a_2=1$ with $b_1,b_2<0$ and $b_1b_2>1$; (iii) $a_1a_2=-1$ with $b_1<0<b_2$ and $b_1b_2<-1$ or with $b_1+b_2<0$ and $b_1-b_2<-2$.)

We shall henceforth be concerned only with the dynamics that has previously been missing, that which occurs inside the discontinuity itself and gives rise to determinacy-breaking at the singularity. In any event, the first step in studying the local dynamics must be to extend the system \eref{2f} in such a way that is can be solved across the discontinuity $x_1=0$, as follows.

%\begin{figure}[h!]\centering%\includegraphics[width=\textwidth]{figso/2folds}
%\vspace{-0.3cm}\caption{\sf Determinacy-breaking at a two-fold singularity of invisible, visible, and bivisible type.}\label{fig:2folds}\end{figure}

%%%%%%%%%%%%%%%%%% %%%%%%%%%%%%%%%%%%
%\subsection{Resolving the discontinuous dynamics}\label{sec:resolve}
%
%
%A theory of structural stability of \eref{gen} relies crucially on an ability to resolve dynamics on the switching threshold $h=0$, by means of re-writing \eref{gen} as a differential inclusion
%\begin{equation}\label{nF}
%{\bf\dot x}={\bf F}\bb{{\bf x};\lambda}\quad:\quad\left\{\begin{array}{lll}\lambda=\sign\bb{h}&\rm if&h\neq0\;,\\\lambda\in\sq{-1,+1}&\rm if&h=0\;.\end{array}\right.
%\end{equation}
%For the local canonical form of the two-fold singularity \eref{2f} this satisfies ${\bf F}\bb{{\bf x};+1}=(-x_2,a_1,b_1)$ and ${\bf F}\bb{{\bf x};-1}=(x_3,b_2,a_2)$. 
%
%All that remains to complete the picture is a set of dynamical equations governing the values of $\lambda$ inside the switching surface. 

%%%%%%%%%%%%%%%%% %%%%%%%%%%%%%%%%%%%%%%%%%%%%%%%%%
\section{Blow up}\label{sec:blowup}

%The procedure for studying the piecewise smooth flow is as follows. 
We shall outline the procedure for studying \eref{2f} in a general form first. Let ${\bf x}=(x_1,x_2,x_3)$ and ${\bf f}=(f_1,f_2,f_3)$, and begin with a piecewise-smooth system 
\begin{equation}\label{fns}
\dot{\bf x}=\left\{\begin{array}{lll}{\bf f}({\bf x};+1)&\rm if&x_1>0\;,\\{\bf f}({\bf x};-1)&\rm if&x_1<0\;,\end{array}\right.
\end{equation}
%and lThe function ${\bf f}({\bf x};\lambda)$ can then be written generally in the form (see \cite{j13error})
which can be extended to $x_1=0$ by defining $\bf f$ as
\begin{equation}\label{fnon}
{\bf f}({\bf x};\lambda)=\frac{1+\lambda}2{\bf f}({\bf x};+1)+\frac{1-\lambda}2{\bf f}({\bf x};-1)+(1-\lambda^2){\bf g}({\bf x};\lambda)
\end{equation}
where
\begin{equation}\label{nF}
\lambda\in\left\{\begin{array}{lll}\sign\bb{x_1}&\rm if&x_1\neq0\;,\\\sq{-1,+1}&\rm if&x_1=0\;.\end{array}\right.
\end{equation}
The function ${\bf g}$ is smooth, and provides ${\bf f}$ with a nonlinear dependence on $\lambda$. Without $\bf g$ the function ${\bf f}({\bf x};\lambda)$ is the convex combination as used in Filippov's standard approach to switching \cite{f88}, and here we shall introduce ${\bf g}\neq0$ only as a small perturbation when necessary. %Below we shall use normal forms for ${\bf f}({\bf x};\pm1)$ at a two-fold singularity. 

The dynamics outside the discontinuity ($x_1\neq0$) is now given by taking either of $\lambda=\pm1$ in \eref{fnon}, which reproduces \eref{fns}. On the discontinuity ($x_1=0$) we consider whether there exists $\lambda\in[-1,+1]$ such that $\dot x_1=0$, for which the flow of \eref{fnon} lies in the tangent plane of the discontinuity, in which case {\it sliding dynamics} along the discontinuity surface is possible. The sliding dynamics can be reached in finite forward or backward time, meaning the flow can lose or gain a dimension as it arrives at or departs from $x_1=0$, so the resulting flow might not be unique. Therefore we shall follow \cite{j13error} and regain uniqueness by {\it blowing up} the discontinuity manifold $x_1=0$, to study the dynamics of $\lambda$ that transports the flow through the jump between $-1$ and $+1$. This is done as follows. 

Because $\lambda$ is a function of $x_1$ only, the dynamics of $\lambda$ is induced by the $x_1$ component of the flow such that $\lambda'=f_1(0,x_2,x_3;\lambda)$, where the prime denotes differentiation with respect to some dummy timescale, instantaneous compared to the main timescale. One way to describe this is to say $\lambda'\equiv\eps\frac{d\lambda}{dt}$ for infinitesimal $\eps>0$. (In fact, this blow up is equivalent to defining $\lambda$ as a function $\lambda=\phi(v)$ of a fast variable $v=x_1/\eps$ such that $\phi(x_1/\eps)\rightarrow{\rm sign}(x_1)$ as $\eps\rightarrow0$; we discuss this further in \sref{sec:conc}). We obtain for the blow up system on $x_1=0$,
\begin{equation}\label{blowup}
(\lambda',\dot x_2,\dot x_3)=\bb{f_1(0,x_2,x_3;\lambda),\;f_2(0,x_2,x_3;\lambda),\;f_3(0,x_2,x_3;\lambda)}\;.
\end{equation}
The steady states of the `fast' prime system satisfy a differential-algebraic system
\begin{equation}\label{blowup0}
(0,\dot x_2,\dot x_3)=\bb{f_1(0,x_2,x_3;\lambda),\;f_2(0,x_2,x_3;\lambda),\;f_3(0,x_2,x_3;\lambda)}\;,
\end{equation}
which defines dynamics on an invariant manifold $$\op M^S=\cc{(\lambda,x_2,x_3)\in[-1,+1]\times\mathbb R^2:\;f_1(0,x_2,x_3;\lambda)=0}\;.$$
The system \eref{blowup0} thus defines {\it sliding dynamics} in the discontinuity surface $x_1=0$, fixing the value of $\lambda\in[-1,+1]$ inside the manifold $\op M^S$ (when it exists) for which $x_1=0$ has invariant dynamics, and specifying the variation $(\dot x_2,\dot x_3)$ on that manifold. In regions of $x_1=0$ where $\op M^S$ does not exist (where $f_1\neq0$ for $x_1=0$ and $\lambda\in[-1,+1]$), the system \eref{blowup} conveys the flow across the discontinuity from one region $x_1\gtrless0$ to the other. 

The systems \eref{fnon} and \eref{blowup} are sufficient to specify the local dynamics, and we shall now apply these to the normal form. The final step in our analysis will be to transform the system \eref{blowup} into the normal form of a folded singularity, but we shall find that this is only possible if ${\bf g}$ is nonzero, and hence there is a nonlinear dependence on $\lambda$.

%%%%%%%%%%%%%%%%% %%%%%%%%%%%%%%%%%%%%%%%%%%%%%%%%%
\section{The unperturbed system}\label{sec:2fdummy}

In this section we show the following. 
\begin{proposition}
The blow-up \eref{blowup} of the normal form two-fold singularity \eref{2f} with $\bf g\equiv0$ is structurally unstable. 
\end{proposition}

To prove this we will perform the blow up described in \sref{sec:intro}, and show that the sliding manifold $\op M^S$ exhibits a degeneracy. 

To resolve the dynamics across the discontinuity, we use \eref{fnon} to express the system \eref{2f} as a convex combination
\begin{eqnarray}\label{2fg0}
(\dot x_1,\dot x_2,\dot x_3)&=&\frac{1+\lambda}2(-x_2,a_1,b_1)+\frac{1-\lambda}2(x_3,b_2,a_2)\\
&:=&\bb{f_1(x_1,x_2,x_3;\lambda),\;f_2(x_1,x_2,x_3;\lambda),\;f_3(x_1,x_2,x_3;\lambda)}\;.\nonumber
\end{eqnarray}
%where $\lambda={\rm sign}(x_1)$ for $h\neq0$. 
Using the blow-up \eref{blowup} on $x_1=0$ we obtain the two timescale system
\begin{equation}\label{h0g0}
(\lambda',\dot x_2,\dot x_3)=(f_1,f_2,f_3)=\frac{1+\lambda}2(-x_2,a_1,b_1)+\frac{1-\lambda}2(x_3,b_2,a_2)\;,
\end{equation}
with $\lambda\in[-1,+1]$. The fast system on $\lambda'$ has equilibria where $\lambda'=0$, yielding a differential-algebraic system 
\begin{equation}\label{blowslidingg0}
(0,\dot x_2,\dot x_3)=(f_1,f_2,f_3)=\frac{1+\lambda}2(-x_2,a_1,b_1)+\frac{1-\lambda}2(x_3,b_2,a_2)\;.
\end{equation}
This describes states that evolve inside the switching surface $x_1=0$ on the main timescale, and these are precisely Filippov's {\it sliding modes}. We can solve in the sliding mode to get $\lambda=(x_3-x_2)/(x_3+x_2)$, and sliding modes exist only when this lies in the allowed range $\lambda\in[-1,+1]$, hence for $x_2x_3>0$. 

In the absence of sliding modes, where $x_2x_3<0$, equation \eref{h0g0} describes the instantaneous transition from one boundary of $\lambda\in\sq{-1,+1}$ to another, whereby the flow crosses through the switching surface as $\lambda$ flips sign.

\begin{figure}[h!]\centering\includegraphics[width=\textwidth]{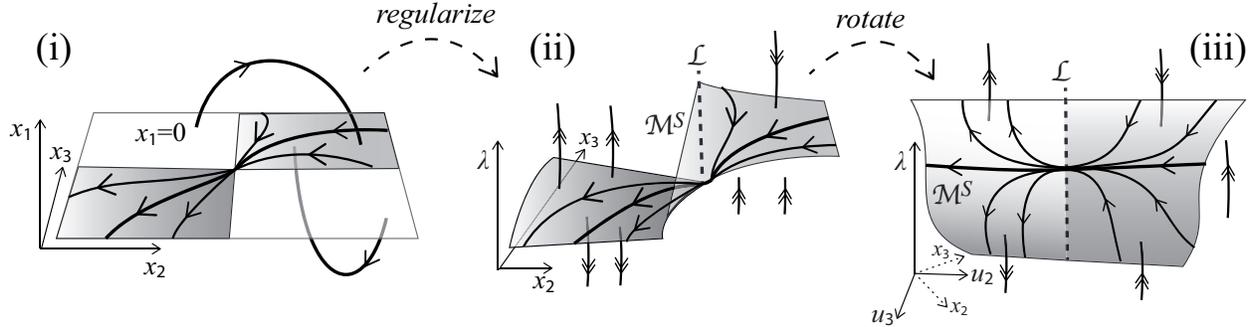}
\vspace{-0.3cm}\caption{\small\sf Blowing up the unperturbed system \eref{2fg0}, for the example of an invisible two-fold. (i) The flow directions outside $x_1=0$  create an attracting sliding region in $x_2,x_3>0$ and repelling sliding region in $x_2,x_3<0$. (ii) The blow up of $x_1=0$ into the interval $\lambda_1\in[-1,+1]$, where the sliding regions become an invariant manifold $\op M^S$ (shaded), hyperbolic except along the vertical line $\op L$, which aligns with the fast (double arrowed) $\lambda_1$ dynamics. (iii) The dynamics in the manifold is best viewed along the $u_3$ axis of rotated coordinates $u_2=x_2+x_3$, $u_3=x_2-x_3$. }\label{fig:2fu}\end{figure}

Concerning the sliding regions $x_2x_3>0$ on $x_1=0$, the subsystem \eref{blowslidingg0} inhabits invariant manifolds of the blow-up system \eref{h0g0} given by
\begin{equation}\label{M2f}
\op M^S=\cc{\;(\lambda,x_2,x_3)\in\sq{-1,+1}\times\mathbb R^2:\;\lambda=\frac{x_3-x_2}{x_3+x_2},\;x_1=0<x_2x_3\;}\;.
\end{equation}
We call $\op M^S$ the {\it sliding manifold}. It consists of two normally hyperbolic branches, one attractive in $x_2,x_3>0$ since $\partial\lambda'/\partial\lambda=-(x_3+x_2)/2<0$, and one repulsive in $x_2,x_3<0$ since $\partial\lambda'/\partial\lambda=-(x_3+x_2)/2>0$. 
The two branches are connected at $x_2=x_3=0$, along a line
\begin{equation}\label{L}
\op L=\cc{(\lambda,x_2,x_3)\in\sq{-1,+1}\times\mathbb R^2:\;x_2=x_3=0,\;}\;.
\end{equation}
Along $\op L$ we have $\partial\lambda'/\partial\lambda=0$, so the manifold $\op M^S$ is no longer normally hyperbolic where it meets $\op L$, and therefore is not invariant there. \Fref{fig:2fu} shows an example of the discontinuous system (i), and its blow-up (ii) showing $\op M^S$ and $\op L$, rotated in (iii) to show $\op L$ more clearly. 

The line $\op L$ represents the two-fold singularity ($x_1=x_2=x_3=0$), stretched out over the interval $\lambda\in[-1,+1]$ in the blow up. Comparing the alignment of $\op L$ with the dynamics outside $\op M^S$, however, reveals a degeneracy in the system \eref{h0g0}. The non-hyperbolic set $\op L$ is a line with tangent vector $e_{\op L}=(1,0,0)$ in the space of $(\lambda,x_2,x_3)$, which means it lies everywhere parallel to the fast one-dimensional $\lambda'$ system in \eref{h0g0}. As a result, all along $\op L$ we have not only $\partial f_1/\partial\lambda=0$, but also $\partial^2f_1/\partial\lambda^2=0$ and moreover $\partial^rf_1/\partial\lambda^r=0$ for all $r>1$. This constitutes an infinite codimension degeneracy, so the system is structurally unstable. This is the result in the proposition.

\bigskip
\noindent{\bf Remarks on the degeneracy}

The degeneracy means that the flow is not transverse to the set $\op L$ where the two sheets of $\op M^S$ intersect. As a result, if we perturb the system, transversal intersection between the attracting and repelling branches of the invariant manifold $\op M^S$ cannot be guaranteed, see e.g. \cite{w05}. Small changes in the expression \eref{2fg0} may result in a system where the two sheets of $\op M^S$ intersect at points, along intervals, or where no intersections occur at all, dependent on the kind of perturbation. These intersections are particularly important because they are they points where {\it canard} trajectories can flow from the attracting to repelling branches of $\op M^S$ (or {\it faux-canard} trajectories can flow in the opposite direction). Hence in the system above the existence of canards cannot be guaranteed. 

In the literature on smooth two timescale systems (as we have here in $(\lambda,x_2,x_3)$ space), the connection of attracting and repelling branches of a slow invariant manifold has been well studied, leading to a generic canonical form and requisite non-degeneracy conditions as described in \cite{w05}. In our notation those conditions are the existence of a point on $\op L$ where
\begin{equation}
\begin{array}{ccccccc}
f_1=0\;,&\quad&\frac{\partial f_1}{\partial\lambda}=0\;,&\quad&
\frac{\partial f_1\;\;\;\;}{\partial x_{2,3}}\neq0\;,&\quad&\frac{\partial^2 f_1}{\partial\lambda^2}\neq0\;,\end{array}\;,
\end{equation}
the fourth of which is violated everywhere in the system described above. 

Adding constant terms or powers of $x_2$ and $x_3$ to \eref{2fg0} would only move the set $\op L$ in the $(x_2,x_3)$ plane, not remove its degeneracy, %, therefore the dummy system \eref{h0g0} around $\op L$ will remain structurally unstable. 
and terms depending on $x_1$ would vanish at the discontinuity. 
The only recourse to break the degeneracy, particularly to give $\frac{\partial^2 f_1}{\partial\lambda^2}\neq0$, is therefore to add terms nonlinear in $\lambda$. Anything we add to the function $f_1$ in \eref{2fg0} must still give \eref{2f}, so it must vanish outside the switching surface $x_1=0$, i.e. be a perturbation of the form \eref{fnon}, as we introduce in the next section.

%\begin{figure}[h!]\centering\includegraphics[width=0.9\textwidth]{figso/2fold}
%\vspace{-0.3cm}\caption{\sf Sketch of splitting and nonsplitting manifolds }\label{fig:}\end{figure}

%%%%%%%%%%%%%%%%% %%%%%%%%%%%%%%%%%%%%%%%%%%%%%%%%%
\section{The perturbed system}\label{sec:2fpert}

In this section we show the following. 
\begin{proposition}
The blow-up \eref{blowup} of the normal form two-fold singularity \eref{2f} with ${\bf g}=(\alpha,0,0)$ can be transformed into
\begin{eqnarray*}%\label{wechs}
\eps\dot{ x_1}&=& x_2+ x_1^2+\ord{\eps x_1,\eps x_3, x_1 x_3}\\
\dot{ x}_2&=& b x_3+ c x_1+\ord{ x_3^2, x_1 x_3}\\
\dot{ x}_3&=& a+\ord{ x_3, x_1}
\end{eqnarray*}
provided $\alpha\neq0$ for small $\eps>0$, where $a,b,c,$ are real constants, and provided 
%$b_1-b_2>2$ if $a_1=-a_2=1$ and $b_1-b_2<-2$ if $a_1=-a_2=-1$ (no conditions if $a_1=a_2$). 
the conditions $\hf(b_1-b_2)\le1=a_1=-a_2$ or $\hf(b_1-b_2)\ge-1=a_1=-a_2$ do not hold. 
\end{proposition}
\noindent It turns out that the case excluded by the conditions $\pm\hf(b_1-b_2)\le1=\pm a_1=\mp a_2$ is that in which there are no canards or faux-canards, i.e. no orbits of the sliding flow passing through the singularity. %, see e.g. the classification in \cite{jc12}.

The perturbed system we consider, applying \eref{fnon} to the normal form \eref{2f} with ${\bf g}=(\alpha,0,0)$, is
\begin{eqnarray}\label{2Fp}
\bb{\dot x_1,\dot x_2,\dot x_3}&=&\frac{1+\lambda}2\bb{-x_2,a_1,b_1}+\frac{1-\lambda}2\bb{x_3,b_2,a_2}+(1-\lambda^2)(\alpha,0,0)\\
&:=&\bb{F_1(x_1,x_2,x_3;\lambda),F_2(x_1,x_2,x_3;\lambda),F_3(x_1,x_2,x_3;\lambda)}\;,\nonumber
\end{eqnarray}
where $\alpha$ is a constant. 
We shall show that perturbing $\dot x_1$ with a small term $\alpha$ in this way is necessary and sufficient for structural stability. Perturbing $\dot x_2$ or $\dot x_3$ is neither necessary nor sufficient, therefore we leave them unaltered with $F_{2,3}\equiv f_{2,3}$ (referring to $(f_1,f_2,f_3)$ in \eref{2fg0}). The systems \eref{2f}, \eref{2fg0} and \eref{2Fp} are identical for all $x_1\neq0$, when $\lambda=\pm1$. 

On $x_1=0$ we consider the blow up of \eref{2Fp} as before, now giving
\begin{equation}\label{2foldblow}
(\lambda',\dot x_2,\dot x_3)=\frac{1+\lambda}2\bb{-x_2,a_1,b_1}+\frac{1-\lambda}2\bb{x_3,b_2,a_2}+\alpha(1-\lambda^2)\bb{1,0,0}\;.
\end{equation}
This is an $\alpha$-perturbation of \eref{h0g0}. 
The $\lambda'$ system has equilibria on the set
\begin{equation}\label{2foldMSp}
\op M^S=\cc{(\lambda,x_2,x_3)\in\sq{-1,+1}\times\mathbb R^2\;:\;-\frac{1+\lambda}2x_2+\frac{1-\lambda}2x_3+\alpha(1-\lambda^2)=0}\;.
\end{equation}
The surface described by $\op M^S$ forms an invariant manifold of \eref{2foldblow} at all points where it is normally hyperbolic, i.e. where $\frac{\partial\lambda'}{\partial\lambda}\neq0$, which now holds except on the curve
\begin{equation}\label{Lpert}
\op L=\cc{(\lambda,x_2,x_3):\;\lambda=2\frac{2\alpha+x_3-x_2}{x_3+x_2}=-\frac{x_3+x_2}{4\alpha}\in\sq{-1,+1}}\;.
\end{equation}
The curve $\op L$ has tangent vector
$e_{\op L}%=\frac{\partial\;}{\partial\lambda}\op L(x_2)=\frac{\partial\;}{\partial\lambda}\bb{\lambda\;,\;+\alpha(\lambda-1)^2\;,\;-\alpha(\lambda+1)^2}$$
=\bb{1\;,\;+2\alpha(\lambda-1)\;,\;-2\alpha(\lambda+1)}$, which for all $|\lambda|\le1$ is transverse to the flow if $\alpha\neq0$. In particular $\op L$ is now transverse to the dot and prime timescale subsystems in \eref{2foldblow}, so the degeneracy of the unperturbed system has been broken. We show below that this perturbed $\alpha\neq0$ system is structurally stable by mapping the two timescale blow-up system \eref{2foldblow} onto the generic folded singularity of smooth systems. 

While the non-hyperbolic curve $\op L$ is now in a general position with respect to the flow, there may exist a new singularity along $\op L$, where the flow's projection along the $\lambda$-direction onto the nullcline $f_1=0$ is indeterminate. This is the so-called {\it folded singularity}, found as follows. The set $f_1=0$ is a fixed point of the $\lambda'$ subsystem of \eref{2foldblow} the flow. A trajectory which remains on $f_1=0$ for an interval of time satisfies $0=\dot f_1=\bb{\lambda',\dot x_2,\dot x_3}\cdot\bb{\partial/\partial\lambda,\partial/\partial x_2,\partial/\partial x_3}f_1=\lambda'\frac{\partial f_1}{\partial\lambda}+(f_2,f_3)\cdot\frac{\partial f_1\quad}{\partial(x_2,x_3)}$ which rearranges to $\frac{\partial f_1}{\partial\lambda}\lambda'=-(f_2,f_3)\cdot \frac{\partial f_1\quad}{\partial(x_2,x_3)}$% along with $(\dot x_2,\dot x_3)=(f_2,f_3)$
. Thus $\lambda'$ is indeterminate at points where 
\begin{equation}\label{fsing}
0=f_1=\frac{\partial f_1}{\partial\lambda}=(f_2,f_3)\cdot \frac{\partial f_1\quad}{\partial(x_2,x_3)}\;.
\end{equation}
Denoting the value of $f_i$ at the singularity as $f_{is}$, and solving \eref{fsing} via
\begin{eqnarray}
0&=&(f_{2s},f_{3s})\cdot\bb{-\frac{1+\lambda_s}2,\frac{1-\lambda_s}2}\nonumber\\
&=&\mbox{$\bb{\frac{a_1+b_2}2+\frac{a_1-b_2}2\lambda_s\;,\;\frac{b_1+a_2}2+\frac{b_1-a_2}2\lambda_s}\cdot\bb{-\frac{1+\lambda_s}2,\frac{1-\lambda_s}2}$}%\\
%&=&-\frac14\bb{\bb{a_1-a_2-b_1+b_2}+2\bb{a_1+a_2}\lambda_s+\bb{a_1-a_2+b_1-b_2}\lambda_s^2}
\end{eqnarray}
we find that the folded singularity lies at $(\lambda,x_2,x_3)=(\lambda_s,x_{2s},x_{3s})$ where
\begin{eqnarray}
\lambda_s=
%&:=&
%\frac{-(a_1+a_2)\pm\sqrt{(a_1+a_2)^2-(a_1-a_2-b_1+b_2)(a_1-a_2+b_1-b_2)}}{a_1-a_2+b_1-b_2}\nonumber\\&=&
%\frac{-(a_1+a_2)\pm\sqrt{(a_1+a_2)^2-(a_1-a_2)^2+(b_1-b_2)^2}}{a_1-a_2+b_1-b_2}\nonumber\\&=&
%\frac{-(a_1+a_2)\pm\sqrt{4a_1a_2+(b_1-b_2)^2}}{a_1-a_2+b_1-b_2}\nonumber\\&=&
\frac{-\frac{a_1+a_2}{b_1-b_2}\pm\sqrt{1+\frac{4a_1a_2}{(b_1-b_2)^2}}}{1+\frac{a_1-a_2}{b_1-b_2}}\;,
\qquad x_{2s}=+\alpha(\lambda_s-1)^2\;,\qquad x_{3s}=-\alpha(\lambda_s+1)^2\;.
\end{eqnarray}
Noting that $a_1$ and $a_2$ just takes values $\pm1$, we have:
%$$
%\lambda_s=\left\{\begin{array}{cll}
%-\frac{2}{b_1-b_2}\pm\sqrt{1+\frac{4}{(b_1-b_2)^2}}&\rm if&a_1=a_2=1\\
%\frac{2}{b_1-b_2}\pm\sqrt{1+\frac{4}{(b_1-b_2)^2}}&\rm if&a_1=a_2=-1\\
%%\frac{\pm\sqrt{(b_1-b_2)^2-4}}{b_1-b_2+2}=
%\pm\sqrt{\frac{b_1-b_2-2}{b_1-b_2+2}}&\rm if&a_1=-a_2=1\\
%%\frac{\pm\sqrt{(b_1-b_2)^2-4}}{b_1-b_2-2}=
%\pm\sqrt{\frac{b_1-b_2+2}{b_1-b_2-2}}&\rm if&a_1=-a_2=-1
%\end{array}\right.
%$$
%which implies:
\begin{itemize}
\item in the case $a_1=a_2=1$, $\lambda_s=\frac{-2}{b_1-b_2}\pm\sqrt{1+\frac{4}{(b_1-b_2)^2}}$, so there exists a unique $\lambda_s\in[-1,+1]$ for any $b_1$ and $b_2$ (the positive root for $b_1>b_2$, the negative root for $b_1<b_2$);
\item in the case $a_1=a_2=-1$, $\lambda_s=\frac{2}{b_1-b_2}\pm\sqrt{1+\frac{4}{(b_1-b_2)^2}}$, so there exists a unique $\lambda_s\in[-1,+1]$ for any $b_1$ and $b_2$ (the positive root for $b_1<b_2$, the negative root for $b_1>b_2$);
\item in the case $a_1=-a_2=1$, $\lambda_s=\pm\sqrt{\frac{b_1-b_2-2}{b_1-b_2+2}}$, there exist two points with $\lambda_s\in[-1,+1]$ for $b_1-b_2>2$, and no points otherwise. % so that for $b_1-b_2\le2$ an equivalence cannot be formed;
\item in the case $a_1=-a_2=-1$, $\lambda_s=\pm\sqrt{\frac{b_1-b_2+2}{b_1-b_2-2}}$, there exist two points with $\lambda_s\in[-1,+1]$ for $b_1-b_2<-2$, and no points otherwise. % so that for $b_1-b_2\ge-2$ an equivalence cannot be formed.
\end{itemize}
In the cases where $\lambda_s$ is unique we proceed directly to the steps that follow below. In the cases where $\lambda_s$ can take two values we can proceed with the following analysis about each value, and will obtain different constants in the final local expression, i.e. a different folded singularity corresponding to each $\lambda_s$. In the cases when $\lambda_s$ does not exist, no equivalence can be formed; these are the cases when the the two-fold's sliding portrait is of focal type (see \cite{jc12}) where orbits winds around the two-fold but sliding orbits never enter or leave it. So excluding those cases $a_1=-a_2=1$ with $b_1-b_2\le2$ and $a_1=-a_2=-1$ with $b_1-b_2\ge-2$, we can proceed as follows. 

Taking a valid value for $\lambda_s$, a translation puts the singularity at the origin of the coordinates
\begin{equation}
y_1=\lambda-\lambda_s\;,\qquad y_2=x_2-x_{2s}\;,\qquad y_3=x_3-x_{3s}\;, 
\end{equation}
%$$f_1=%\frac{\bar x_{3s}-\bar x_{2s}}2-\frac{\bar x_{3s}+\bar x_{2s}}2\lambda_s+\alpha(1-\lambda_s^2)
%\frac{x_3-x_2}2   -\frac{x_3+x_2}2\lambda_s   -\bb{\frac{\bar x_{3s}+\bar x_{2s}}2+2\alpha\lambda_s}z   -\bb{\frac{x_3+x_2}2+\alpha y_1}y_1$$
then $f_1$ becomes
\begin{equation}\label{fz}
f_1=-\frac{\lambda_s+1}2y_2    -\frac{\lambda_s-1}2y_3    -\bb{\frac{y_3+y_2}2+\alpha y_1}y_1
%=x\cdot f_{1s,x}+y_1(y_2,y_3)\cdot f_{1s,(y_2,y_3)y_1}+\hfy_1^2f^{,y_1y_1}
\end{equation}
with derivatives evaluated at the singularity
 \begin{equation}
 \begin{array}{rclcrcl}
\frac{\partial f_{1}}{\partial y_1}%=x\cdot f^{s,xz}+zf^{,zz}
&=&-\frac{y_{3}+y_{2}}2-2\alpha y_{1}\;,&\qquad&\frac{\partial f_{1}\quad}{\partial(y_2,y_3)}&=&-\bb{\frac{\lambda_s+1}2,\frac{\lambda_s-1}2}\;,\\\frac{\partial^2 f_{1}}{\partial^2 y_1^2}&=&-2\alpha\;,&\qquad&\frac{\partial^2 f_{1}\qquad}{\partial(y_2,y_3)\partial y_1}&=&-\hf(1,1)\;.\end{array}
\end{equation}
%%$$\op L=\cc{(y_1,y_2,y_3):\; \frac{\lambda_s+1}2y_2  +\frac{\lambda_s-1}2y_3  +\bb{\frac{y_3+y_2}2+\alpha y_1}y_1=\frac{y_3+y_2}2+2\alpha y_1=0,\; |y_1+\lambda_s|<1}$$
%%$$\op L=\cc{(y_1,y_2,y_3):\; \frac{\lambda_s(y_2+y_3)+(y_2-y_3)}2  -\alpha y_1^2=\frac{y_3+y_2}2+2\alpha y_1=0,\; |y_1+\lambda_s|<1}$$
%%$$\op L=\cc{(y_1,y_2,y_3):\; \frac{y_2-y_3}2 -2\alpha\lambda_sy_1  -\alpha y_1^2=\frac{y_3+y_2}2+2\alpha y_1=0,\; |y_1+\lambda_s|<1}$$
Solving the conditions $f_{1}=\frac{\partial f_{1}}{\partial y_1}=0$ we find that $\op L$ can be written in functional form as
$$(y_1,y_2,y_3)=\op L(y_1):=\bb{y_1,-\alpha y_1(2 -2\lambda_s - y_1),\; -\alpha y_1(2+2\lambda_s + y_1)}\;.$$
%$$\op L=\cc{(y_1,y_2,y_3):\;x\cdot f_{1s,(y_2,y_3)}+y_1(y_2,y_3)\cdot f_{1s,(y_2,y_3)y_1}+\hfy_1^2f_{1s,y_1y_1}=x\cdot f_{1s,(y_2,y_3)y_1}+y_1f_{1s,y_1y_1}=0}$$
%which can be solved to give on $\op L$
%$$\frac{y_{3L}}\alpha+ 2(1+\lambda_s)y_{1L} + y_{1L}^2=0\;,\qquad y_2%=-\frac{y_3f_{1s,y_3y_1}+y_1f_{1s,y_1y_1}}{f_{1s,y_2y_1}}
%=-y_3-4\alpha y_1\;.$$
It is more useful to parameterize $\op L$ in terms of $y_3$, which we obtain as
\begin{equation}
\op L:\;\;\bb{\begin{array}{c}y_{1L}(y_3)\\y_{2L}(y_3)\\y_{3L}(y_3)\end{array}}:=\bb{ \begin{array}{c}  -1-\lambda_s+\sqrt{(1+\lambda_s)^2-y_3/\alpha} \\ -y_3-4\alpha (-1-\lambda_s+\sqrt{(1+\lambda_s)^2-y_3/\alpha})\\   y_3  \end{array} }\;,
\end{equation}
in terms of three functions $y_{1L},y_{2L},y_{3L},$ whose derivatives are
\begin{equation}
y_{1L}'(y_3)%=\frac{-1/2\alpha}{\sqrt{(1+\lambda_s)^2-y_3/\alpha}}
=\frac{-1/2\alpha}{1+\lambda_s+y_{1L}(y_3)}\;,\quad y_{2L}'(y_3)%=-1-4\alpha\frac{-1/2\alpha}{\sqrt{(1+\lambda_s)^2-y_3/\alpha}}
=\frac{1-\lambda_s-y_{1L}(y_3)}{1+\lambda_s+y_{1L}(y_3)}\;,\quad y_{3L}'(y_3)=1\;.
\end{equation}
We can then rectify $\op L$ to lie along some $z_3$ axis by defining new coordinates
\begin{equation}
z_1=y_1-y_{1L}(y_3)\;,\qquad z_2=y_2-y_{2L}(y_3)\;,\qquad z_3=y_3\;.
\end{equation}
The former vector field components become
\begin{equation}
\begin{array}{rcl}
f_1&=&-\frac{1+\lambda_s+y_{1L}}2z_2  - \alpha z_1^2-z_2z_1/2\;,\\
f_2&=&\frac{a_1+b_2}2+\frac{a_1-b_2}2\lambda_s+\frac{a_1-b_2}2\bb{z_1+y_{1L}(z_3)}=f_{2s}+\ord{z_3,z_1}\;,\\
f_3&=&\frac{b_1+a_2}2+\frac{b_1-a_2}2\lambda_s+\frac{b_1-a_2}2\bb{z_1+y_{1L}(z_3)}=f_{3s}+\ord{z_3,z_1}\;.\end{array}
\end{equation}
When we now evaluate the derivatives $z_1',\dot z_2,\dot z_2,$ it is necessary to relate the dot and prime time derivatives. We do this using an infinitesimal $\eps>0$, letting a dot denote the derivative with respect to $t$, while a prime denotes the derivative with respect to a `fast' time $t/\eps$ which becomes instantaneous as $\eps\rightarrow0$. Then
\begin{eqnarray*}
z_1'&=&%(y_1-y_{1L}(z_3))'=
f_1-\eps\dot z_3y_{1L}'(z_3)\\
%&=&-\frac{y_{1L}(z_3)+\lambda_s+1}2z_2  - \alpha z_1^2-z_2z_1/2+\eps\frac{ f_{3s}+\ord{z_3,z_1}}{2\alpha(y_{1L}(z_3)+1+\lambda_s)}\\
%&=&-\frac{\lambda_s+1}2\bb{z_2-\frac{\eps f_{3s}}{\alpha(1+\lambda_s)^2}}  - \alpha z_1^2 -\frac{y_{1L}(z_3)}2z_2-z_2z_1/2+\eps\ord{z_3,z_1}\\
&=&-\frac{\lambda_s+1}2\bb{z_2-\frac{\eps f_{3s}}{\alpha(1+\lambda_s)^2}}  - \alpha z_1^2+\ord{\eps z_3,\eps z_1,z_2z_3,z_2z_1}\;,
\end{eqnarray*}
using $y_{1L}(z_3)=\ord{z_3}%=-1-\lambda_s+(1+\lambda_s)\sqrt{1-\frac{y_3}{\alpha(1+\lambda_s)}}
=-\frac{z_3}{2\alpha}+\ord{z_3^2}\;$. 
We then make a small corrective shift $\tilde z_2=z_2-\frac{\eps f_{3s}}{\alpha(1+\lambda_s)^2}$ which yields
\begin{eqnarray*}
\dot{\tilde z}_2
&=&%\dot{(y_2-y_{2L}(z_3))}=
f_2-\dot z_3y_{2L}'(z_3)%\\
%&=&f_2- f_3\frac{1-\lambda_s-y_{1L}(y_3)}{1+\lambda_s+y_{1L}(y_3)}\\
%&=&f_{2s}+(z_1+y_{1L}(z_3))f_{2s,y_1}-\sq{f_{3s}+(z_1+y_{1L}(z_3))f_{3s,y_1}}\frac{1-\lambda_s-y_{1L}(y_3)}{1+\lambda_s+y_{1L}(y_3)}\\
%&=&\frac{\sq{f_{2s}+(z_1+y_{1L}(z_3))f_{2s,y_1}}(1+\lambda_s+y_{1L}(y_3))-\sq{f_{3s}+(z_1+y_{1L}(z_3))f_{3s,y_1}}(1-\lambda_s-y_{1L}(y_3))}{1+\lambda_s+y_{1L}(y_3)}\\
%&=&%\frac{(1+\lambda_s)f_{2s}-(1-\lambda_s)f_{3s}}{1+\lambda_s+y_{1L}(y_3)}+
%y_{1L}(y_3)\frac{f_{2s}+f_{3s}}{1+\lambda_s+y_{1L}(y_3)}
%+(z_1+y_{1L}(z_3))\frac{(1+\lambda_s)f_{2s,y_1}-(1-\lambda_s)f_{3s,y_1}}{1+\lambda_s+y_{1L}(y_3)}\\
%&&+y_{1L}(y_3)(z_1+y_{1L}(z_3))\frac{f_{2s,y_1}+f_{3s,y_1}}{1+\lambda_s+y_{1L}(y_3)}\\
%&=&-\frac{z_3}{2\alpha}\frac{f_{2s}+f_{3s}}{1+\lambda_s}
%+(z_1-\frac{z_3}{2\alpha})\frac{(1+\lambda_s)f_{2s,y_1}-(1-\lambda_s)f_{3s,y_1}}{1+\lambda_s}+\ord{z_1z_3,z_3^2}\\
%&=&-\frac{z_3}{2\alpha}\frac{f_{2s}+f_{3s}}{1+\lambda_s}
%-\frac{z_3}{2\alpha}c
%+cz_1+\ord{z_1z_3,z_3^2}\\
=\frac b\alpha z_3+cz_1
+\ord{z_3^2,z_3z_1}
\end{eqnarray*}
where
$$c= \frac{\partial f_{2s}}{\partial\lambda}- \frac{1-\lambda_s}{1+\lambda_s}\frac{\partial f_{3s}}{\partial\lambda}\;,\qquad b=-\frac{1}{2}\bb{\frac{ f_{2s}+ f_{3s}}{1+\lambda_s}+c}\;.$$

The last thing to do is just scaling. Collecting everything together so far we have
\begin{equation}
\begin{array}{rcl}
z_1'&=&d_1{\tilde z}_2-\alpha z_1^2+\ord{\eps z,\eps z_3,z_1z_3}\\
\dot {\tilde z}_2&=&\frac b\alpha z_3+cz_1+\ord{z_3^2,z_1z_3}\\
\dot z_3&=& f_{3s}+\ord{z_3,z_1}
\end{array}
\end{equation}
where $d_1=-\hf(1+\lambda_s)$. 
%and group terms, denoting the derivative with respect to $-t{\rm sign}(\alpha)\sqrt{|\alpha|}/\eps$ with a dagger (so we have to scale all time by $-{\rm sign}(\alpha)$), 
%\begin{eqnarray*}
%(-{\rm sign}(\alpha)d_1\dot {\tilde z}_2)&=&\frac{-d_1b}{|\alpha|} \bb{-{\rm sign}(\alpha)z_3}+\frac{d_1c}{\sqrt{|\alpha|}}\bb{\sqrt{|\alpha|}z_1}+\ord{z_3^2,z_1z_3}\\
%(-{\rm sign}(\alpha)\dot z_3)&=& f_{3s}+\ord{{\tilde z}_2,z_3,z_1}\\
%\sqrt{|\alpha|}z_1^\dagger&=&\bb{-{\rm sign}(\alpha)d_1{\tilde z}_2}+(\sqrt{|\alpha|}z_1)^2+\ord{\eps z,\eps z_3,z_1z_3}
%\end{eqnarray*}
Defining new variables $\tilde x_1=\sqrt{|\alpha|}z_1$, $\tilde x_2=-{\rm sign}(\alpha)d_1 {\tilde z}_2$, $\tilde x_3=-{\rm sign}(\alpha)z_3$, $\tilde t=-{\rm sign}(\alpha)t$, gives
\begin{equation}\label{wechs}
\begin{array}{rcl}
{\tilde x_1}'&=&\tilde x_2+\tilde x_1^2+\ord{\eps x_1,\eps\tilde x_3,\tilde x_1\tilde x_3}\\
\dot{\tilde x}_2&=&\tilde b\tilde x_3+\tilde c\tilde x_1+\ord{\tilde x_3^2,\tilde x_1\tilde x_3}\\
\dot{\tilde x}_3&=&\tilde a+\ord{\tilde x_3,\tilde x_1}
\end{array}
\end{equation}
where
%$$c= \frac{\partial f_{2s}}{\partial\lambda}- \frac{1-\lambda_s}{1+\lambda_s}\frac{\partial f_{3s}}{\partial\lambda}\;,\qquad b=-\frac{1}{2}\bb{\frac{ f_{2s}+ f_{3s}}{1+\lambda_s}+c}\;,\qquad d_1=-\hf(1+\lambda_s)\;,$$ 
%$$\tilde a= f_{3s}\;,\qquad\tilde b=-\frac{d_1b}{|\alpha|}$$
\begin{equation}\label{wp}
\begin{array}{l}
\tilde a= f_{3s}%=\frac{b_1+a_2}2+\frac{b_1-a_2}2\lambda_s
\;,\quad\tilde b%=-\frac{d_1b}{|\alpha|}%=-\frac{1}{4|\alpha|}(1+\lambda_s)\bb{\frac{ f_{2s}+ f_{3s}}{1+\lambda_s}+c}
%=-\frac{1}{4|\alpha|}\bb{ f_{2s}+ f_{3s}+{(\lambda_s+1) \frac{\partial f_{2s}}{\partial\lambda}+(\lambda_s-1) \frac{\partial f_{3s}}{\partial\lambda}}}
=-\frac{1}{4|\alpha|}\bb{f_{2s}+ f_{3s}-2\tilde c\sqrt{|\alpha|}}\;,\\
\tilde c%=\frac{d_1c}{\sqrt{|\alpha|}}
=\frac{-1}{2\sqrt{|\alpha|}}\bb{(\lambda_s+1)\frac{\partial f_{2s}}{\partial\lambda}+(\lambda_s-1)\frac{\partial f_{3s}}{\partial\lambda}}
\;.
\end{array}
\end{equation}
Since $x'\equiv\eps\dot x$ this is the result in the proposition, and it is clearly valid only for $\alpha\neq0$.

\begin{figure}[h!]\centering\includegraphics[width=\textwidth]{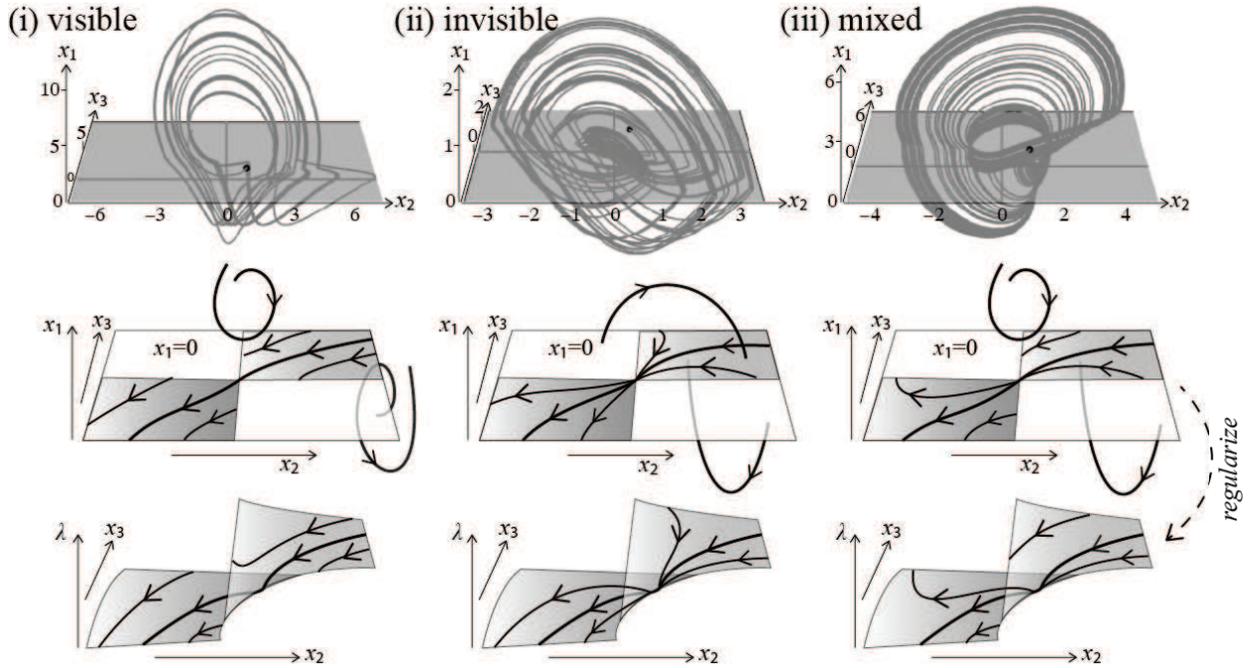}
\vspace{-0.3cm}\caption{\small\sf Three examples of attractor organised around a two-fold singularity. Showing: (i) a simulated trajectory, (ii) a sketch of the discontinuous flow inside and outside $x_1=0$ that gives rise to it, and (iii) the blow up on $x_1=0$.}\label{fig:2fs}\end{figure}

\Fref{fig:2fs} shows an example of the perturbed system and its blow up for each flavour of two-fold in (i) (corresponding to those in \fref{fig:2folddet}), followed by their blow up (ii), and a rotation (iii) to show the phase portrait around the set $\op L$ more clearly (similar to \fref{fig:2fu}). 
In the most extreme case, the folded node, the original phase portrait contains infinitely many intersecting trajectories traversing the singularity, while the perturbed system splits these into distinguishable orbits, a finite number of which asymptote to the attracting and repelling branches of the critical manifold.  

\begin{figure}[h!]\centering\includegraphics[width=\textwidth]{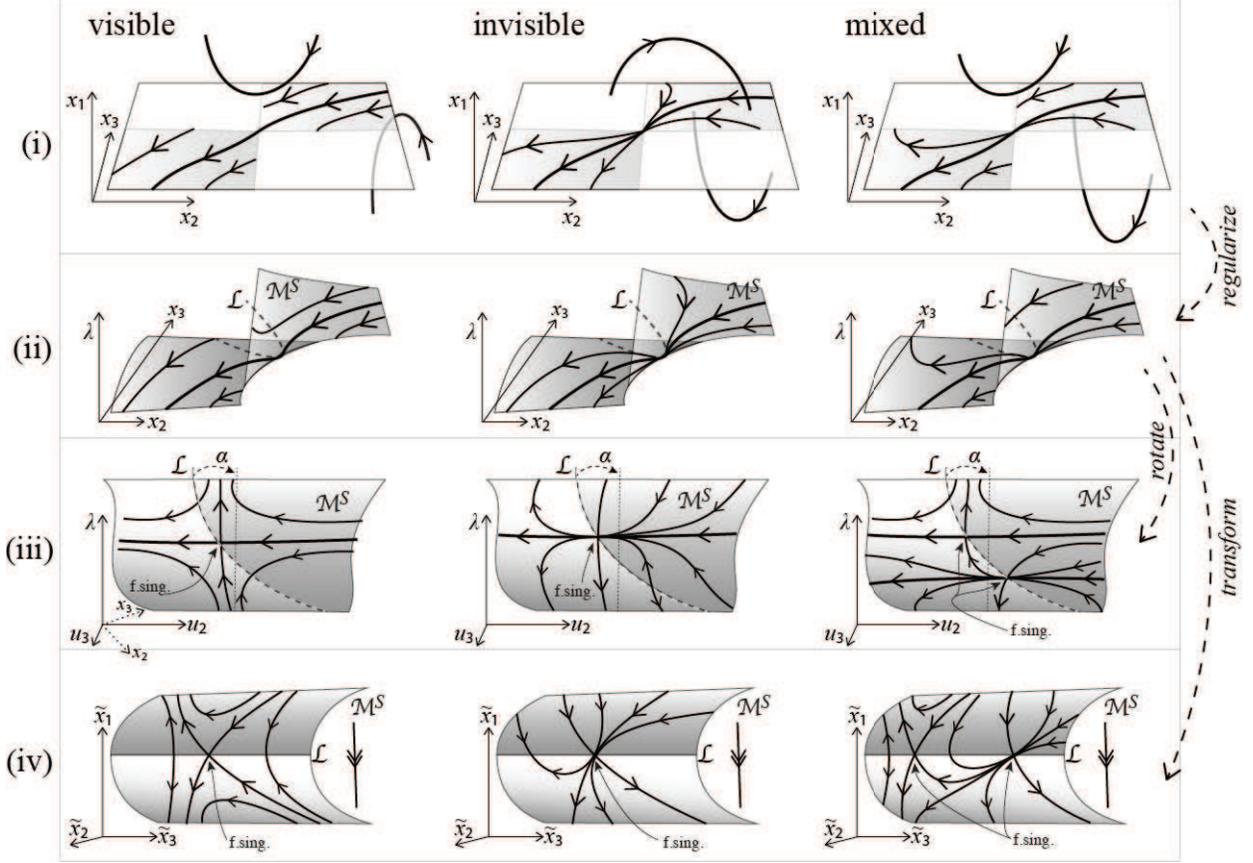}
\vspace{-0.3cm}\caption{\small\sf Blowing up the perturbed ($\alpha\neq0$) system, for examples of each flavour of two-fold. Labelling as in \fref{fig:2fu}. Note in the blow up (ii) that $\op L$ is now a curve. Rotating around the $\lambda$ axis in (iii) we can see the attracting branch (upper right segment) and repelling branch (lower left segment) of the invariant manifold $\op M^S$ (shaded), connected by $\op L$. The folded singularity (f.sing.) appears along $\op L$, two in the case of mixed visibility, recognised as having a phase portrait that resembles a saddle or node if we reverse time in the repelling branch of $\op M^S$. In (iv) we sketch the corresponding phase portraits in the slow-fast system \eref{wechs}. 
%(i) The flow directions outside $x_1=0$  create an attracting sliding region in $x_2,x_3>0$ and repelling sliding region in $x_2,x_3<0$. (ii) The blow up of $x_1=0$ into the interval $\lambda_1\in[-1,+1]$, where the sliding regions become an invariant manifold $\op M^S$ (shaded), hyperbolic except along the vertical line $\op L$, which aligns with the fast (double arrowed) $\lambda_1$ dynamics. (iii) The dynamics in the manifold is best viewed along the $u_3$ axis of rotated coordinates $u_2=x_2+x_3$, $u_3=x_2-x_3$.
}\label{fig:2fr}\end{figure}

\bigskip
\noindent{\bf Remarks on the singularity}

A glance at the papers \cite{w05,w12,des10} reveals what a charismatic singularity lies hidden in the dynamics of the two-fold, revealed by blow-up, and waiting to be released when the discontinuous system is perturbed by simulations that smooth, regularize, or otherwise approximate the discontinuity. 

In particular one may ask what happens in the cases of determinacy-breaking illustrated in \fref{fig:2folddet}, since a smooth flow should be deterministic. Because the two-fold maps to a folded singularity (actually one or two folded singularities as we saw above), if we consider what happens with $\eps>0$ we find determinacy is restored, but replaced by sensitivity to initial conditions. 

The manifold $\op M^S$ becomes the slow critical manifold of a slow-fast system, and for $\eps$ by standard results of geometric singular perturbation theory \cite{f79}, there exist invariant manifolds $\op M^S_\eps$ in the $\eps$-neighbourhood of $\op M^S$ but away from $\op L$, on which the dynamics is topologically equivalent to the sliding dynamics (termed `reduced' or slow dynamics usually in this context) found above. Trajectories that pass close to the singularity, or more precisely, close to the folded singularity on the set $\op L$, may persist in following the  manifold $\op M^S$ from its stable to unstable branches, while other nearby trajectories will veer wildly away, their fate sensitive to initial conditions and proximity to primary canard orbits (those which persist along both branches of $\op M^S$ throughout the local region). 

%$$A=\f^{20}/2\alpha$$

%$$1>\lambda_s^2=\frac{1+\frac{4a_1a_2}{(b_1-b_2)^2}+\bb{\frac{a_1+a_2}{b_1-b_2}}^2\mp2\frac{a_1+a_2}{b_1-b_2}\sqrt{1+\frac{4a_1a_2}{(b_1-b_2)^2}}}{1+2\frac{a_1-a_2}{b_1-b_2}+\bb{\frac{a_1-a_2}{b_1-b_2}}^2}$$
%$$\bb{b_1-b_2}\bb{a_1-a_2}-4a_1a_2>\mp\bb{a_1+a_2}\sqrt{(b_1-b_2)^2+4a_1a_2}$$
%$$\left\{\begin{array}{l}
%2<\pm\sqrt{(b_1-b_2)^2+4}\\
%2<\mp\sqrt{(b_1-b_2)^2+4}\\
%b_2-b_1<2\\
%b_1-b_2<2
%\end{array}\right.\qquad\Rightarrow\qquad\left\{\begin{array}{l}
%0<(b_1-b_2)^2\\
%0<(b_1-b_2)^2\\
%b_2-b_1<2\\
%b_1-b_2<2
%\end{array}\right.$$
%
%for real $\lambda_s$
%%$$1+\frac{4a_1a_2}{(b_1-b_2)^2}>0$$
%$$(b_1-b_2)^2+4a_1a_2>0$$
%
%%$$1+\lambda_s=\frac{1-2\frac{a_2}{b_1-b_2}\pm\sqrt{1+\frac{4a_1a_2}{(b_1-b_2)^2}}}{1+\frac{a_1-a_2}{b_1-b_2}}$$
%
%

Like the different kinds of two-fold, there are different classes of folded singularity, and their classification depends on the slow dynamics inside $\op M^S$. From the expressions \eref{wechs}-\eref{wp} we see that the class therefore depends not only on the constants $a_1,a_2,b_1,b_2,$ of the original piecewise smooth system, but also on the `hidden' parameter $\alpha$.  
%
%$$\op M^S=\cc{(\tilde x_1,\tilde x_2,\tilde x_3)\in\times\mathbb R^3:\;\tilde x_2=-\tilde x_1^2,\;\abs{\frac{\tilde x_1}{\sqrt{|\alpha|}}+\lambda_s+y_{1L}(-{\rm sign}(\alpha)z_3)}\le1}$$

The classification scheme is fairly simple, and can be used to verify the dynamics on $\op M^S$ seen in \fref{fig:2fr}. The projection of the system \eref{wechs} onto $\op M^S$, found by differentiating the condition $0=\tilde x_2+\tilde x_1^2$ with respect to time to give $0=\tilde b\tilde x_3+\tilde c\tilde x_1+2\tilde x_1\dot{\tilde x}_1+\ord{\tilde x_3^2,\tilde x_1\tilde x_3}$, is
$$\bb{\begin{array}{c}\dot{\tilde x}_1\\\dot{\tilde x}_3\end{array}}=\frac1{-2\tilde x_1}\cc{\bb{\begin{array}{cc}\tilde c&\tilde b\\-2\tilde a&0\end{array}}\bb{\begin{array}{c}\tilde x_1\\\tilde x_3\end{array}}+\ord{\tilde x_3^2,\tilde x_1\tilde x_3}}\;.$$
A classification then follows by neglecting the singular prefactor $1/2\tilde x_1$ and considering whether the phase portrait is that of a focus, a node, or a saddle. This is determined by the $2\times2$ matrix Jacobian, which has 
trace $\tilde c$, determinant $2\tilde a\tilde b$, and eigenvalues $\hf(\tilde c\pm\sqrt{\tilde c^2-8\tilde a\tilde b})$. 
This will not be the true system's phase portrait because the time-scaling from the $1/2\tilde x_1$ factor is positive in the attractive branch of $\op M^S$, negative (time-reversing) in the repulsive branch, and divergent at the singularity (turning infinite time convergence to the singularity into finite time passage {\it through} the singularity). The effect of this is to `fold' together attracting and repelling pairs of each equilibrium type, so each equilibrium becomes a `folded-equilibrium', forming a continuous bridge between branches of $\op M^S$. 
As a result the flow on $\op M^S$ is a folded-saddle if $\tilde a\tilde b<0$, a folded-node if $0<8\tilde a\tilde b<\tilde c^2$, and a folded-focus if $\tilde c^2<8\tilde a\tilde b$. Canard cases occur for $\tilde c>0$ and faux canard for $\tilde c<0$. In \fref{fig:2fr} we show the result of blowing up the discontinuous system for an example of each type of two-fold that exhibits determinacy-breaking (those from \fref{fig:2folddet}). In the visible two-fold the singularity becomes a folded-saddle, in the invisible case it becomes a folded-node, while the mixed case becomes a pair consisting of one folded-saddle and one folded-node.

One may ask why certain cases, in particular those without determinacy-breaking, were excluded by the proposition above. These are systems where the sliding phase portraits contain no orbits passing between the attracting and repelling sliding regions, and therefore one expects that in their blow up there should exists no canards connecting attracting and repelling branches of $\op M^S$, and hence no folded singularities. If there are no such connections then no intersections between the two branches of $\op M^S$ persist when we let $\eps$ be nonzero. Hence the omission of these cases is consistent, and a posteori it is obviously necessary. %The result, that the blow up of a discontinuous two-fold where no orbits passing between the attracting and repelling sliding regions, results in a smooth system with non-intersecting attracting and repelling slow manifolds, could be proven more fully, but fits intuitively with the picture now emerging. 

\begin{figure}[h!]\centering\includegraphics[width=\textwidth]{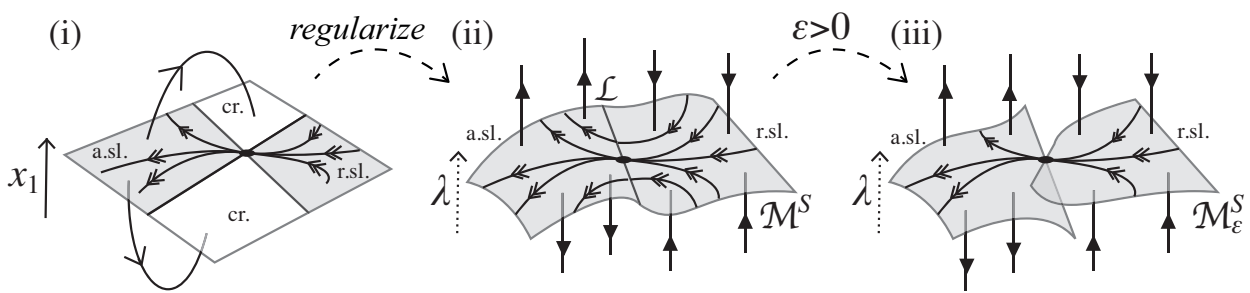}
\vspace{-0.3cm}\caption{\small\sf Sketch of: (i) a two-fold singularity; (ii) its blow up with $\alpha\neq0$; and (iii) the splitting of the branches of $\op M^S$ when we replace $ x_1'$ with $\eps\dot x_1$ and consider small $\eps>0$. (iii) approaches (ii) as $\eps\rightarrow0$.}\label{fig:degen}\end{figure}

In \fref{fig:degen} we illustrate the way splitting of manifolds occurs when we introduce the instantaneous timescale (whose derivative is denoted with a prime) with a fast timescale $t/\eps$ for $\eps>0$. The case shown is an invisible two-fold. 

An in-depth description of the dynamics that ensues in the different cases of two-folds, and the smoothings subject to perturbations, would be lengthy, and deserves future study elsewhere. 
As a demonstration we conclude with three examples of complex oscillatory attractors formed by two-fold singularities, simulated by smoothing. We take $\alpha=1/5$ and:
\begin{enumerate}
\item[(i)] ${\bf f}^+=(-x_2,\sfrac25x_1+\sfrac1{10}x_2-1,\sfrac3{10}x_2-\sfrac15x_2x_3-\sfrac25)$, ${\bf f}^-=(x_3,\sfrac15x_2x_3-\sfrac35,\sfrac25x_3-1-x_1)$;
\item[(ii)] ${\bf f}^+=(-x_2,1+x_1,-\sfrac75)$, ${\bf f}^-=(x_3,-\sfrac9{10},1-\sfrac35x_1)$;
\item[(iii)] ${\bf f}^+=(-x_2+\sfrac1{10}x_1,x_1-\sfrac65,x_1-2)$, ${\bf f}^-=(x_3+\sfrac1{10}x_1,x_1+\sfrac{23}{100},1-x_1)$;
\end{enumerate}
simulated in \fref{fig:2fs}. 
The simulation uses Mathematica's {\sf NDSolve}, for which the discontinuous function $\lambda_1={\rm sign}\;x_1$ is replaced with a sigmoid function $\tanh(x_1/\eps)$ with $\eps=10^{-5}$. Because the system is structurally stable, the exact choice of smoothing function is not important so long as it is simple, for example continuous and monotonic, and further simulations not shown here verify that such dynamics persists with different smoothing functions, for example varying $\eps$, or replacing $x_1$ with a different sigmoid function such as $x_1/\sqrt{\eps^2+x_1^2}$. The first row shows the simulation of a single trajectory for a time interval $t=1000$ in $(x_1,x_2,x_3)$ space, as the flow attempts to switch between the critical slow manifold flow in (ii) given by the blow up on $x_1=0$, and the $x_1>0$ and $x_1<0$ flows in (i). In the second row we plot the familiar picture of sliding dynamics on $x_1=0$ for the unperturbed system (neglecting nonlinear $\lambda$), and sketch the curvature of the $x_1\neq0$ flows. The sliding vector field has canard trajectories passing from the righthand attracting branch to the lefthand repelling branch via the folded singularity; at a visible two-fold only one canard exists, at an invisible two-fold every sliding trajectory is a canard, and at a mixed two-fold a region of trajectories are canards. The third row shows the critical manifold $\op M^S$ in the $(\lambda_1,x_2,x_3)$ blow up.

%%%%%%%%%%%%%%%%% %%%%%%%%%%%%%%%%%%
%%%%%%%%%%%%%%%%% %%%%%%%%%%%%%%%%%%
\section{Closing Remarks}\label{sec:conc}

For such a simple system, even taking its piecewise linear local normal form, the two-fold exhibits intricate and varied dynamics. Just how intricate becomes even more clear as we attempt to regularize the discontinuity, and study how the two-fold related to slow-fast dynamics of smooth systems. As well as insight into the dynamics that is seen upon simulated such a system, % have seen, attempting to resolve dynamics at the singularity $x_1=x_2=x_3=0$ involves even more intricacy, and 
this adds a new facet to the question of the structural stability of the two-fold, which has remained a stimulating question since \cite{t90}. 

The degeneracy in \sref{sec:2fdummy} gives some insight into why simulations of systems containing two-fold singularities with determinacy-breaking singularities are highly unpredictable. In attempting to solve near a discontinuity, simulations may introduce aspects of smoothing, hysteresis, delay and noise, not represented in the ideal discontinuous model, with unpredictable effects on such a structurally unstable high codimension degeneracy. We therefore have to resolve this degeneracy before expecting to obtain comprehensible simulations, and without leaving the class of piecewise smooth systems, the nonlinear switching term introduced in the perturbed system in \sref{sec:2fpert} is sufficient. 

The system \eref{2Fp} now succeeds \eref{2fg0} as our prototype for the local dynamics of the two-fold singularity (with $\lambda$ given by \eref{nF} in both cases). The question of whether this constitutes a `normal form' has issues both in the discontinuous and slow-fast settings, but it is clear that \eref{2Fp} is structurally stable, and represents all classes of behaviour that occur both in the discontinuous system, and in its blow up to a slow-fast system. 

Finally we should emphasize that the results of this paper are equivalent if we smooth, rather than blow up, the discontinuity. In particular the degeneracy that makes the perturbation ${\bf g}\neq0$ necessary is %not considered in previous papers %(e.g. \cite{teix/hogan}) which employ smoothing directly rather than blow up, but the degeneracy is 
still present. To relate the procedure of smoothing (or Sotomayor-Teixeira's regularization \cite{ts11}) to the blow up procedure, we consider $\lambda$ as a function rather than a dynamic variable. To smooth \eref{2fg0} we let $\lambda=\phi(x_1/\eps)$ for small $\eps>0$, then define a new variable $v=\eps x_1$ for which $v'=\eps\dot v=f_1(\eps v,x_2,x_3;\phi(v))$, which replaces $\lambda'$ in \eref{h0g0} to give us our two timescale system. Note that $\lambda'$ and $v'$ are related by scalar, $\lambda'=\eps\dot\lambda=\eps\phi'(v)\dot v$, amounting only to a time rescaling. Hence the $\phi$-smoothed system is qualitatively similar to the $\lambda$-blow up system, with the same invariant manifold $\op M^S$, and containing the same line $\op L$ which lies aong the fast $v'$ flow when $\eps=0$.

%\begin{figure}[h!]\centering\includegraphics[width=\textwidth]{figso/./figures/}
%\vspace{-0.3cm}\caption{\sf  }\label{fig:}\end{figure}

\bibliography{../grazcat}

\begin{thebibliography}{10}

\bibitem{b81}
E.~Beno\^{i}t, J.~L. Callot, F.~Diener, and M.~Diener.
\newblock Chasse au canard.
\newblock {\em Collect. Math.}, 31-32:37--119, 1981.

\bibitem{cj10}
A.~Colombo and M.~R. Jeffrey.
\newblock Non-deterministic chaos, and the two-fold singularity in piecewise
  smooth flows.
\newblock {\em SIAM J. App. Dyn. Sys.}, 10:423--451, 2011.

\bibitem{jc12}
A.~Colombo and M.~R. Jeffrey.
\newblock The two-fold singularity: leading order dynamics in n-dimensions.
\newblock {\em Physica D}, 265:1--10, 2013.

\bibitem{jd10}
M.~Desroches and M.~R. Jeffrey.
\newblock Canards and curvature: nonsmooth approximation by pinching.
\newblock {\em Nonlinearity}, 24:1655--1682, 2011.

\bibitem{jd12}
M.~Desroches and M.~R. Jeffrey.
\newblock Nonsmooth analogues of slow-fast dynamics: pinching at a folded node.
\newblock {\em submitted}, 2013.

\bibitem{des10}
M.~Desroches, B.~Krauskopf, and H.~M. Osinga.
\newblock Numerical continuation of canard orbits in slow-fast dynamical
  systems.
\newblock {\em Nonlinearity}, 23(3):739--765, 2010.

\bibitem{f79}
N.~Fenichel.
\newblock Geometric singular perturbation theory.
\newblock {\em J. Differ. Equ.}, 31:53--98, 1979.

\bibitem{f88}
A.~F. Filippov.
\newblock {\em Differential Equations with Discontinuous Righthand Sides}.
\newblock Kluwer Academic Publ. Dortrecht, 1988.

\bibitem{j13error}
M.~R. Jeffrey.
\newblock Hidden dynamics in models of discontinuity and switching.
\newblock {\em Physica D}, 273-274:34--45, 2014.

\bibitem{jc09}
M.~R. Jeffrey and A.~Colombo.
\newblock The two-fold singularity of discontinuous vector fields.
\newblock {\em SIAM Journal on Applied Dynamical Systems}, 8(2):624--640, 2009.

\bibitem{t90}
M.~A. Teixeira.
\newblock Stability conditions for discontinuous vector fields.
\newblock {\em J. Differ. Equ.}, 88:15--29, 1990.

\bibitem{t93}
M.~A. Teixeira.
\newblock Generic bifurcation of sliding vector fields.
\newblock {\em J.Math.Anal.Appl.}, 176:436--457, 1993.

\bibitem{ts11}
M.~A. Teixeira and P.~R. da~Silva.
\newblock Regularization and singular perturbation techniques for non-smooth
  systems.
\newblock {\em Physica D}, 241(22):1948--55, 2012.

\bibitem{w05}
M.~Wechselberger.
\newblock Existence and bifurcation of canards in $\mathbb{R}^3$ in the case of
  a folded node.
\newblock {\em SIAM J. App. Dyn. Sys.}, 4(1):101--139, 2005.

\bibitem{w12}
M.~Wechselberger.
\newblock A propos de canards (apropos canards).
\newblock {\em Trans. Amer. Math. Soc}, 364:3289--3309, 2012.

\end{thebibliography}
\bibliographystyle{plain}

\clearpage

\end{document}